\documentclass{article}
\usepackage{arxiv}
\usepackage{amsmath}
\usepackage{cleveref}
\usepackage{comment}
\usepackage{tikz}

\usepackage[english]{babel}
\usepackage[utf8x]{inputenc}
\usepackage[T1]{fontenc}

\usepackage{mathrsfs}
\usepackage{graphicx}

\usepackage{multirow}

\usepackage{tikz}
\usetikzlibrary{shapes}
\usetikzlibrary{patterns}
\usetikzlibrary{arrows,positioning,decorations.pathmorphing,trees}

\usepackage{url}

\usepackage{chngcntr}
\counterwithin{table}{section}
\counterwithin{figure}{section}

\usepackage{amsfonts}

\usepackage[ruled,vlined]{algorithm2e}

\usepackage{enumerate}

\usepackage{caption}
\usepackage{subcaption}

\providecommand{\keywords}[1]
{
  \small	
  \textbf{\textit{Keywords---}} #1
}

\newcommand{\weights}{weights}
\newcommand{\weight}{weight}

\title{Penalties and Rewards for Fair Learning in Paired Kidney Exchange Programs}

\author{Margarida Carvalho\thanks{CIRRELT \& D\'epartement d'informatique et de recherche op\'erationnelle, Universit\'e de Montr\'eal: carvalho@iro.umontreal.ca}
\and {\bf Alison Caulfield}\thanks{Department of Mathematics \& Statistics, McGill University: alison.caulfield@mail.mcgill.ca}
\and {\bf Yi Lin}\thanks{Department of Mathematics \& Statistics, McGill University: yi.lin2@mail.mcgill.ca}
\and {\bf Adrian Vetta}\thanks{School of Computer Science and Dept. of Mathematics \& Statistics, McGill University: adrian.vetta@mcgill.ca}}

\begin{document}

\maketitle
\begin{abstract}
A kidney exchange program, also called a kidney paired donation program, can be viewed as a repeated, dynamic trading and 
allocation mechanism. This suggests that a dynamic algorithm for transplant exchange selection
may have superior performance in comparison to the repeated use of a static algorithm.
We confirm this hypothesis using a full scale simulation of the Canadian Kidney Paired 
Donation Program: learning algorithms, that attempt to learn optimal patient-donor 
 {\weights} in advance via dynamic simulations, do lead to improved outcomes. 
Specifically, our learning algorithms, designed with the objective of 
fairness (that is, equity in terms of transplant accessibility across cPRA groups), also lead to an increased number of transplants 
and shorter average waiting times. Indeed, our highest performing learning algorithm improves egalitarian 
fairness by 10\% whilst also increasing the number of transplants by 6\% 
and decreasing waiting times by 24\%. However, our main result is much more surprising.
We find that the most critical factor in determining the performance 
of a kidney exchange program is 
{\bf not} the judicious assignment of positive {\weights} (rewards) to patient-donor pairs.
Rather, the key factor in increasing the number of transplants, decreasing waiting times 
and improving group fairness is the judicious assignment of a negative {\weight} (penalty) to the 
small number of non-directed donors in the kidney exchange program.
\end{abstract}


\section{Introduction}

The Canadian Kidney Paired Donation (KPD) program is a kidney exchange program (KEP) consisting of incompatible patient-donor pairs and non-directed anonymous donors, henceforth called \emph{nodes}. The program uses a weighting (point) scheme to select compatible kidney exchanges between the donors and patients~\cite{CBS18}. Our work was motivated by the accumulation over time of \emph{hard-to-match} patients in the KPD program.
We ask whether or not this phenomenon, a recurring problem in KEPs~\cite{DPS14}, can be mitigated via the use of dynamic learning algorithms to select arc/node {\weights} in the kidney exchange graphs used for transplant selection.
We answer this question in the affirmative. Our key finding is that a judicious choice of {\weights} 
can simultaneously improve performance with respect to multiple criteria including group fairness,
waiting times, and the number of transplants. However, this conclusion requires  
that negative {\weights} (penalties) are permitted in transplant selection, in addition to positive {\weights} (rewards).

\subsection{Background and Related Literature}\label{sec:background}
Kidney exchanges were first proposed by Rapaport~\cite{Rap86} and the first programs 
were instigated in South Korea in the 1990s~\cite{PMK99,HKJ08}.
Recently, these programs have spread across the globe: see~\cite{MM14} and~\cite{GKW14} for details on 
the national kidney exchange programs in the UK and the Netherlands, respectively; see also~\cite{BHA19} for an 
overview of KEPs in Europe. Canada started its own KEP, 
entitled the \emph{Kidney Paired Donation (KPD) Program}, in 2008~\cite{CNC15}.

\begin{figure}[h!]
     \centering
     \begin{subfigure}[b]{0.34\textwidth}
         \centering
         \begin{tikzpicture}[scale=.45]
        \draw (0,4) circle [radius=.75]; \node at (0,4) {{\tiny $p_1,d_1$}};
        \draw (4,4) circle [radius=.75]; \node at (4,4) {{\tiny $p_2,d_2$}};
        \draw (0,0) circle [radius=.75]; \node at (0,0) {{\tiny $p_3,d_3$}};
        \draw   [fill=gray!20] (4,0) circle [radius=.75]; \node at (4,0) {{\tiny $\emptyset,d_{4}$}};
        \draw [<-] (.75,0) -- (3.25,0);
        \draw [ ->] (.55,.55) -- (3.45,3.45);
        \draw [<-] (0,.75) -- (0,3.25); 
        \draw [<-] (.75,4) -- (3.25, 4); 
        \end{tikzpicture}
        \centering
        \caption{A kidney exchange graph.}\label{fig:kepgraph}
     \end{subfigure}
     \hfill
     \begin{subfigure}[b]{0.28\textwidth}
         \centering
         \begin{tikzpicture}[scale=.45]
            \draw (0,4) circle [radius=.75]; \node at (0,4) {{\tiny $p_1,d_1$}};
            \draw (4,4) circle [radius=.75]; \node at (4,4) {{\tiny $p_2,d_2$}};
            \draw (0,0) circle [radius=.75]; \node at (0,0) {{\tiny $p_3,d_3$}};
            \draw   [fill=gray!20] (4,0) circle [radius=.75]; \node at (4,0) {{\tiny $\emptyset,d_{4}$}};
            \draw [<-] (.75,0) -- (3.25,0);
            \draw [ultra thick, ->] (.55,.55) -- (3.45,3.45);
            \draw [ultra thick, <-] (0,.75) -- (0,3.25); 
            \draw [ultra thick, <-] (.75,4) -- (3.25, 4); 
            \end{tikzpicture}
            \centering
            \caption{A transplant cycle.}\label{fig:cycle}
     \end{subfigure}
     \hfill
     \begin{subfigure}[b]{0.34\textwidth}
         \centering
         \begin{tikzpicture}[scale=.45]
        \draw (0,4) circle [radius=.75]; \node at (0,4) {{\tiny $p_1,d_1$}};
        \draw (4,4) circle [radius=.75]; \node at (4,4) {{\tiny $p_2,d_2$}};
        \draw (0,0) circle [radius=.75]; \node at (0,0) {{\tiny $p_3,d_3$}};
        \draw   [fill=gray!20] (4,0) circle [radius=.75]; \node at (4,0) {{\tiny $\emptyset,d_{4}$}};
        \draw [ultra thick, <-] (.75,0) -- (3.25,0);
        \draw [ultra thick, ->] (.55,.55) -- (3.45,3.45);
        \draw [<-] (0,.75) -- (0,3.25); 
        \draw [ultra thick, <-] (.75,4) -- (3.25, 4); 
        \end{tikzpicture}
        \centering
        \caption{A transplant path (chain).}\label{fig:path}
     \end{subfigure}
        \caption{Illustration of feasible exchanges for a KEP.}
        \label{fig:three graphs}
\end{figure}
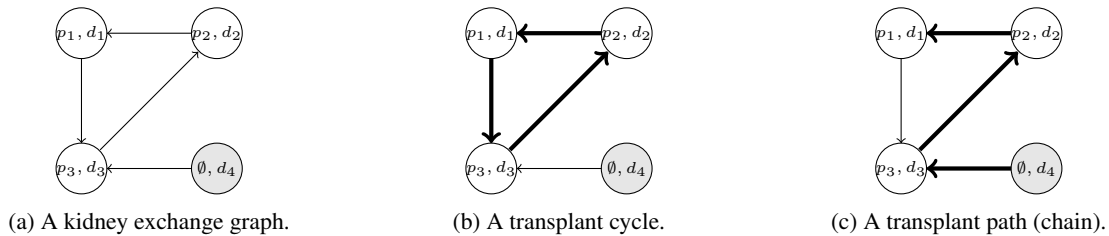

The theoretical foundations underlying kidney exchange programs were provided by Roth et al.~\cite{RSU04}; see 
also S\"onmez and \"Unver~\cite{SU13}. Succinctly, a KEP can be modelled by a directed graph, called a {\em kidney exchange graph} 
whose nodes consist of incompatible patient-donor pairs, denoted $(p_i, d_i)$, and non-directed anonymous donors (NDADs), denoted
$(\emptyset, d_j)$, also called altruistic donors. 
An arc in the graph indicates that a transplant is feasible between the corresponding donor and patient.
We refer the reader to Section~\ref{sec:model} for a detailed description of kidney exchange graphs.
For now, we remark that kidney exchanges are selected using either a directed cycle or a directed simple path in the graph.
 For an illustration consider Figure~\ref{fig:three graphs}, which shows a KEP with
three incompatible patient donor pairs and one non-directed anonymous donor (shaded).
A valid directed cycle consists only of patient-donor pairs; see Figure~\ref{fig:cycle}.
A valid directed path must contain exactly one non-directed anonymous donor located at the source of the path; see Figure~\ref{fig:path}. 
Notice that when the cycle is selected, the altruistic donor $d_4$ is unused and 
remains available for future transplants. On the other hand, with the path, the donor $d_1$ is unused.
The importance of directed paths and, thus, of altruistic donors in transplant selection 
was originally highlighted by Montgomery et al.~\cite{MGM06} and Roth et al.~\cite{RSU06}. 
Moreover, a detailed investigation into the most efficacious use of altruistic donors will be a key
focus of this work. 

The big question is what transplants (arcs) should be selected from the graph to optimize health care outcomes.
At first thought, this appears obvious: maximize the number of transplants. This is the approach taken in KEPs that use 
the myopic algorithm for transplant selection. There are three major reasons why it is questionable. First, transplant selection is a multi-criteria decision problem, of which the quantity of transplants is just one aspect.
For example, medical practitioners must also take into account the 
health of the patient and the urgency of the case, waiting times (both individual and collective), 
the quality of potential transplants, etc. Thus, choosing to maximize the number of transplants above other objectives is a subjective decision. Second, the myopic algorithm is implicitly unfair. It will prioritize patients who are easy-to-match (e.g., 
in terms of blood-group compatibility and antibody compatibility) at the expense of hard-to-match patients.
Third, KEPs are repeated, dynamic mechanisms. This, in turn, has several consequences.
Perhaps counter-intuitively, algorithms to optimize an objective in the short-term may 
not optimize the objective over the long-term. A striking example of this is the use of static pricing mechanisms 
in electricity markets. There, pricing for short-term efficiency can lead to huge inefficiencies in the long term~\cite{CM10}.
For KEPs, this means the myopic approach of maximizing the number of transplants in each round may
not maximize the number of transplants over many years. Our experiments will show this fact clearly.
The dynamic nature of KEPs can also have deleterious effects. For example, a serious problem currently 
facing KEPs is the accumulation over time of hard-to-match patients. 

A standard approach, and the one taken in the KPD program~\cite{CNC15}, to attempt to address these issues is to add weights (also called points or rewards) to prioritize certain types of exchanges.
For instance, higher weight may be assigned to hard-to-match patients such as those who are highly-sensitized\footnote{Patients with high calculated Panel Reactive Antibody (cPRA) rates.} or have 
O-blood type. Dickerson et al.~\cite{DPS12} proposed the use of weights and in~\cite{DPS14}, they explain the usefulness of weights in dynamic contexts and in dealing with different objectives; see also~\cite{DS15,FBS20}. 
Given a weighting, the resultant optimization problem is to find a
maximum weight collection of node-disjoint cycles and paths in the kidney exchange graph. 
The first integer programming (IP) formulations were introduced by Abraham et al.~\cite{ABT07} and by Roth et al.~\cite{RSU07}, and the first compact IP was presented in~\cite{CON13}.
The optimization problem is NP-hard, but state-of-the-art software works well in practice. 
For example, the position-index formulations 
proposed in~\cite{DMP16} can be efficiently solved by out-of-the-box IP solvers when path lengths are capped, while the TSP-based approach of~\cite{AAG15} is preferred when path lengths are uncapped. When only cycles are considered, the most efficient algorithm in practice is the decomposition approach by~\cite{LAM20}. 
 In addition to arc/node weightings, hierarchical optimization techniques have been used to incorporate prioritization criteria~\cite{DGG20}; see Biro et al.~\cite{BHA19,BKM21} for a overview of weighted and hierarchical approaches in European KEPs. 

So weights are important in KEPs, but what is the best choice of weights? 
The myopic approach corresponds to assigning equal weight to each patient-donor node. An alternative approach is
expert-selected weights, where medical practitioners formulate a set of rules by which   
patient-donor nodes and/or arcs are prioritized (weighted). 
We apply a third approach: given the dynamic nature of a KEP, we use dynamic learning algorithms to 
select the weights. But, as alluded to, for multi-criteria decision problems such as 
in KEPs, the ``best'' choice of weights will depend upon the somewhat subjective choice of what criteria
to optimize. Our focus in this paper is to learn weights to optimize {\em group fairness}.
Informally, we will group the patients according to how hard they are to match and attempt to design
transplant selection algorithms that are fair for each group (see Section~\ref{sec:fair} for details).
In particular, we wish to reduce the accumulation over time of hard-to-match patients. 
Fortunately, as will be shown, group fairness can be achieved whilst also balancing the objectives of 
maximizing transplant numbers and minimizing waiting times.

\subsection{Overview and Results}\label{sec:overview}
In Section~\ref{sec:model} we describe kidney exchange programs and explain how to model them using
a kidney exchange graph. 
In Section~\ref{sec:algs} we present a collection of learning algorithms 
designed to learn weights with the aim of optimizing group fairness.
The weights are learnt using a full-scale simulation of end-stage kidney disease patients in Canada.
We also use the simulator to assess the performance of these algorithms in comparison
to both the myopic algorithm and the Canadian KPD algorithm.
For both tasks, learning the weights and performance testing, we use the implementation of the position-index formulations
provided in~\cite{DMP16}. 

In Section~\ref{sec:fair} we present a specific class of group fairness measures which 
will be used to analyze the performance of 
our transplant selection algorithms.
The three most important measures in the class are {\em Nash group fairness}, {\em utilitarian group fairness}, and {\em egalitarian group fairness}; indeed, the latter measure is particularly suited for the task of evaluating equity in transplant accessibility.

The results from our experiments are described in Section~\ref{sec:results}.
All of our learning algorithms significantly outperform the myopic algorithm and the Canadian KPD algorithm; see Table~\ref{tab:main}.
In particular, our overall highest performing algorithm improves the egalitarian fairness measure by 10\% over the Canadian KPD algorithm
and by 21\% over the myopic algorithm. That is perhaps unsurprising given that the algorithms were designed to optimize fairness.
However, the learning algorithm also leads to an increase of 6\% in the number of transplants, a 
decrease of 24\% in waiting times, and a 2\% reduction in the use of altruistic donors. This is much more surprising. In particular, the 
myopic algorithm was designed to maximize the number of transplants but the learning algorithms do significant better in
the long run for that criterion. This illustrates the need to incorporate the repeated and dynamic nature of KEPs into 
transplant selection algorithms. 

As stated, our experiments were tested using Canadian population statistics, but
we believe the take-home lessons from these experiments apply to the design of KEPs more generally.
These include the following six lessons:\\[-0.3cm]

\noindent {\tt Lesson I.} {\em There is negligible benefit in allowing for paths of length greater than 5.} 
\hfill [Table~\ref{tab:paths}] \\[-0.3cm]

\noindent {\tt Lesson II.} {\em It is very important to allow for cycles of length 3. The improvements with cycles of lengths 4 or 5 are marginal.}\hfill [Table~\ref{tab:cycles}] \\[-0.3cm]

\noindent {\tt Lesson III.} {\em Assigning an appropriate negative weight (penalty) to each altruistic donor is fundamental in  determining the quality of outcomes.} \hfill [Table~\ref{tab:altruist-weights} and Table~\ref{tab:prices}]\\[-0.3cm]

\noindent {\tt Lesson IV.} {\em The marginal contribution in terms of lives saved per altruistic donor is non-increasing
in the number of altruists. (At the current altruist arrival rate in Canada, each altruistic donor saves approximately 
two lives.)} \hfill [Table~\ref{tab:altruist-number}] \\[-0.3cm]

\noindent {\tt Lesson V.} {\em Fair algorithms can be used to reduce the accumulation of hard-to-match patients in waiting pools.}  \hfill [Table~\ref{tab:fairness}]\\[-0.3cm]

\noindent {\tt Lesson VI.} {\em Fair algorithms can be implemented without detrimental affects on utilitarian criteria such as
waiting times and the number of transplants.} \hfill [Table~\ref{tab:main} and Table~\ref{tab:prices}]\\[-0.3cm]

We remark the first two lessons will not surprise practitioners. 
Whilst the potential importance of long paths was originally highlighted by Ashlagi et al.~\cite{AGRR12}, 
early works gave theoretical and experimental evidence that consideration of only short cycles~\cite{RSU07} or short paths~\cite{DPS12b} suffice.

The last two lessons validate the approach taken in this paper. Lesson V confirms that learning algorithms can 
reduce the accumulation of hard-to-match patients. Moreover improved fairness need not come at a cost; our fair 
learning algorithms can be used to simultaneously improve the performance of a KEP over a range a measures. 

What is surprising, at least to us, is that Lesson III turns out to be by far the most important of the lessons.
Put simply, weighting (negatively) altruistic donors has much greater beneficial effects than weighting (positively) all the other nodes combined! To our knowledge this fact has not previously been discovered. 
Indeed, as we show in Table~\ref{tab:prices}, simply incorporating a negative weight or penalty for each altruistic donor in
the myopic and KPD algorithms produces outcomes that are comparable to our learning algorithms.
We remark that this lesson contrasts with recent work~\cite{DGHR18} giving theoretical arguments for the
prioritization of paths over cycles.

Lesson IV is useful because it provides for an understanding of how to exploit Lesson III for maximum benefit.
Furthermore, there is evidence~\cite{DPS12b} that some altruistic donors prefer a shorter active time on a deceased donor 
program list rather than a longer active time on a living donor program list. Knowledge that remaining in the KEP will save approximately two lives rather than just one in the deceased donor program may help alleviate this issue.\footnote{Of course, if an altruistic donor requests only a short wait before donation then that can easily be accommodated in our model with the use of individual weight adjustments. However, such requests 
induce serious ethical questions~\cite{PCP11,WDA10}.}

\section{The Graphical Kidney Exchange Model}\label{sec:model}

We now present a dynamic model of a kidney exchange.
First, we describe in Section~\ref{sec:static} 
a graphical model with a static set of patients and donors. In particular, we explain how feasible transplants (representing compatible donors and patients) can be modelled by an exchange graph.
Then, in Section~\ref{sec:transplant-selection}, we summarise the main objectives driving transplant selection in this static model. Section~\ref{sec:dynamic} concludes with a dynamic model that mimics the evolution over time of the patient and donor pool in the kidney exchange program.

\subsection{A Static Model}\label{sec:static}
A {\em Kidney Exchange Program} (KEP) can be modelled by a directed graph $G = (V, A)$, called an {\em exchange graph}.
A node $i\in V$ is either a {\em patient-donor pair} or is a {\em non-directed anonymous donor} (NDAD), also called an {\em altruistic donor}.
In the former case, $i=\{p_i, d_i\}$ consists of a patient $p_i$ and an accompanying (incompatible) donor $d_i$.
In the later case, $i=\{\emptyset, d_i\}$ consists simply of a donor with no accompanying patient.
There is an arc $(i,j) \in A$ if the kidney of donor $d_i$ is compatible with patient~$p_j$.
Note that a node~$i$ corresponding to an altruistic donor must have in-degree zero in $G$.

It follows that a directed cycle in $G$ of patient-donor pairs forms a feasible exchange of donors and thus, a feasible set of transplants.
Furthermore, a directed path, also called a {\em chain}, whose source is an altruistic donor also forms 
a feasible set of transplants. Therefore, a feasible set of transplants is a node-disjoint collection
of paths and cycles.
An example of an exchange graph is shown in Figure~\ref{fig:exchange-graph}, where a feasible
set of transplants is shown in bold.

\begin{figure}[h!]
\centering
\begin{minipage}[t]{0.45\textwidth}
\begin{tikzpicture}[scale=0.4]
\draw  (0,15) circle [radius=1.35]; \node at (0,15) {{\footnotesize $p_1,d_1$}};
\draw [fill=gray!20]  (15,10) circle [radius=1.35]; \node at (15,10) {{\footnotesize $\emptyset,d_8$}};
\draw (10,5) circle [radius=1.35]; \node at (10,5) {{\footnotesize $p_{11},d_{11}$}};

\draw [<-] (1.35,15) -- (3.65,15); 
\draw [->] (0,11.35) -- (0,13.65); 
\draw [->] (6.35,15) -- (8.65,15); 
\draw [->] (5,13.65) -- (5,11.35); 
\draw [->] (15,11.35) -- (15,13.65); 
\draw [->] (15,8.65) -- (15,6.35); 
\draw [->] (15,3.65) -- (15,1.35); 
\draw [->] (1,1) -- (4,4);
\draw [->] (6,6) -- (9,9);
\draw [->] (6,4) -- (9,1);
\draw [->] (11,4) -- (14,1);
\draw [->] (9.5,6.25) -- (9.5,8.75);
\draw [<-] (10.5,6.25) -- (10.5,8.75);
\draw [<-] (11.35,10) -- (13.65,10);

\draw  [fill=gray!20] (0,0) circle [radius=1.35]; \node at (0,0) {{\footnotesize $\emptyset,d_{13}$}};
\draw (5,0) circle [radius=1.35]; \node at (5,0) {{\footnotesize $p_{14},d_{14}$}};
\draw (10,0) circle [radius=1.35]; \node at (10,0) {{\footnotesize $p_{15},d_{15}$}};
\draw (15,0) circle [radius=1.35]; \node at (15,0) {{\footnotesize $p_{16},d_{16}$}};
\draw [ultra thick, ->] (1.35,0) -- (3.65,0); 
\draw [ultra thick, ->] (6.35,0) -- (8.65,0); 
\draw [ultra thick, ->] (11.35,0) -- (13.65,0); 

\draw  [fill=gray!20] (5,15) circle [radius=1.35]; \node at (5,15) {{\footnotesize $\emptyset,d_{2}$}};
\draw (10,10) circle [radius=1.35]; \node at (10,10) {{\footnotesize $p_7,d_7$}};
\draw (15,5) circle [radius=1.35]; \node at (15,5) {{\footnotesize $p_{12},d_{12}$}};
\draw [ultra thick, ->] (6,14) -- (9,11); 
\draw [ultra thick, ->] (11,9) -- (14,6); 

\draw   (10,15) circle [radius=1.35]; \node at (10,15) {{\footnotesize $p_{3},d_{3}$}};
\draw  (15,15) circle [radius=1.35]; \node at (15,15) {{\footnotesize $p_{4},d_{4}$}};
\draw [ultra thick, <-] (11.25,14.5) -- (13.75,14.5); 
\draw [ultra thick, ->] (11.25,15.5) -- (13.75,15.5); 

\draw (0,5) circle [radius=1.35]; \node at (0,5) {{\footnotesize $p_9,d_9$}};
\draw (5,5) circle [radius=1.35]; \node at (5,5) {{\footnotesize $p_{10},d_{10}$}};
\draw (0,10) circle [radius=1.35]; \node at (0,10) {{\footnotesize $p_5,d_5$}};
\draw (5,10) circle [radius=1.35]; \node at (5,10) {{\footnotesize $p_6,d_6$}};
\draw [ultra thick, <-] (1.35,5) -- (3.65, 5); 
\draw [ultra thick, ->] (1.35,10) -- (3.65, 10); 
\draw [ultra thick, ->] (0, 6.35) -- (0, 8.65); 
\draw [ultra thick, <-] (5, 6.35) -- (5, 8.65);
\end{tikzpicture}
\caption{An exchange graph and a feasible set of transplants.}\label{fig:exchange-graph}
\end{minipage}
\centering
\begin{minipage}[t]{0.45\textwidth}
\begin{center}
\begin{tikzpicture}[scale=0.4]
\draw  (5,10) circle [radius=1]; \node at (5,10) {O};
\draw  (5,0) circle [radius=1]; \node at (5,0) {AB};
\draw  (0,5) circle [radius=1]; \node at (0,5) {A};
\draw  (10,5) circle [radius=1]; \node at (10,5) {B};

\draw [thick, ->] (5,9) -- (5,1.2); 
\draw [thick, ->]  (5.7,9.3) -- (9.7,6.2); 
\draw [thick, ->] (4.3,9.3) -- (0,6.2); 
\draw [thick, ->] (0,4) -- (4.4,1); 
\draw [thick, ->] (10,4) -- (5.7,1); 
\end{tikzpicture}
\end{center}
\caption{Blood-group compatibilities.}\label{fig:blood-types}
\end{minipage}
\end{figure}

It will be important to understand the factors that determine transplant compatibility and, hence, the arcs in the exchange graph.
The first major factor is {\em blood-group compatibility}. 
There are four blood groups, namely $O,A,B$ and $AB$.
A donor of group $O$ is a {\em universal donor} and may donate to a patient of any blood group.
A donor of group~$A$ may donate to a patient of group $A$ or $AB$.
A donor of group~$B$ may donate to a patient of group $B$ or $AB$.
Finally, a group $AB$ may donate only to a patient also of group $AB$.
This is illustrated in Figure~\ref{fig:blood-types}, where the arrows point from donor to feasible recipients of other blood groups.
From the perspective of a patient this means that a patient of group $AB$ is a {\em universal recipient} as it may 
receive a transplant from a patient of any blood group. In contrast, a patient of group $O$ can only receive a transplant
from a donor also of group $O$. 
In terms of the exchange graph, a patient-donor pair $\{p_i, d_i\}$ of blood groups $\{AB, O\}$, respectively,
can, in principle, have an outgoing arc to every other patient-donor pair and an incoming arc from 
every other patient-donor pair. However, a patient-donor pair $\{p_i, d_i\}$ of blood groups $\{O, AB\}$, respectively,
can only have an outgoing arc to a node whose patient is of blood group $AB$ and can only have 
an incoming arc from a node whose donor is of blood group $O$.
The second major factor concerns {\em antibodies} and {\em antigens}. 
Antibodies protect against foreign pathogens.
However, if an antibody in the patient recognises and binds to an antigen from the donor this can lead to graft failure. 
The larger the number of types of antibodies in the patient
the less likely a donor will be compatible with that patient. 
This is measured by the cPRA rate which estimates the proportion of donors against which the patient has antibodies rendering the match incompatible: easy-to-match patients have a cPRA close to zero and
hard-to-match patients close to one.

\subsection{Transplant Selection}\label{sec:transplant-selection}
 The general approach  for transplant selection is to add
 a weight $w_i$ to each node $i$ or a weight $w_{ij}$ to each arc $(i,j)$ in the graph.\footnote{It is easy to see that arc weights are more general than node weights.
Given a node weight, we can implement the same system using arc weights as follows.
If we want node $j$ to have weight $w_j$ then simply assign $w_{ij}=w_j$ for each arc $(i,j)$.} 
 This induces an optimization problem whose optimal solution is a set of 
node-disjoint cycles and paths of maximum total weight. In addition, it is natural to impose a maximum cycle length, $C$,
and a maximum path length, $P$, where the length of a cycle/path is the corresponding number of patients. 
Three basic methods are used in selecting the weights. 
First, the simplest approach uses $w_{i}=1$ for every patient-donor pair $i$ (and $w_{i}=0$ for every altruistic donor $i$).
The resultant objective is then to maximize the total number of transplants
subject to the cycle and path length constraints; the corresponding algorithms are called {\em myopic}.
The second approach uses {\em expert-defined} weights.
Here the arc/node weights are prescribed in advance by medical experts to signify the relative 
priority and importance of a given transplant. 
This is the method implemented by the Canadian KPD program (Section~\ref{sec:algs}).
The third approach uses learned-weights. The weights are learnt internally by the system 
to optimize some performance objective~\cite{DPS12}.
This will be the approach taken in this paper where our objective will be group fairness. 

\subsection{A Dynamic Model}\label{sec:dynamic}
We can formulate and solve this kidney-exchange model. But this model is static, that is, it only encompasses
one-time period. 
This is problematic, since non-selected patients are still waiting for transplants, as is readily seen from Figure~\ref{fig:exchange-graph}, where the patient-donor pairs  $\{p_1, d_1\}$ and $\{p_{11}, d_{11}\}$ are still waiting for transplants 
and the altruistic donor $\{\emptyset, d_8\}$ is still willing to donate a kidney.
Moreover, over time there are new node arrivals to the program 
and also unforeseen departures from the  program (e.g. patients who die or receive a kidney from an alternate source, 
paired-donors who drop out, etc).  
Consequently, in reality, a KEP is dynamic with exchanges calculated at regular time intervals.
For example, in Canada, matches are calculated every $4$ months. In particular, in time period $t$ there is an
exchange graph $G_t=(V_t, A_t)$. In time period $t+1$ there is a new exchange graph $G_{t+1}=(V_{t+1}, A_{t+1})$ 
including new arrivals but excluding nodes matched in the previous round and other departures. 

To build a dynamic kidney-exchange model we need, in particular, arrival and departure rates for each type of patient-donor pairs and altruistic donors. 
To calculate these we used data from~\cite{CBS18}.
In particular, the expected arrival rate of patient-donor pairs per time-period (4 months) is $37$ which
we model with a random Poisson distribution with mean $37$. We also model the arrival rate of altruistic donors
through an independent Poisson process with mean $\lambda_A=4.5625$.
The blood type probabilities of donors and patients are
$\mathbb{P}(A, B, AB, O) \ =\  (0.46, 0.42, 0.09, 0.03).$ 
Furthermore, the probabilities that the cPRA rate of a patient lie in specific intervals are
$\mathbb{P}\big((0, 0), (0.01, 0.50), (0.51, 0.94), (0.95, 0.96), (0.97, 1)\big)=(0.24, 0.29, 0.24, 0.10, 0.13)$.
Denote these interval probabilities by $(\alpha_1, \alpha_2, \alpha_3, \alpha_4, \alpha_5)=(0.24, 0.29, 0.24, 0.10, 0.13)$;
these are used in~\cite{CL22} and correspond to the distribution among patients registered in the KPD program when it was launched~\cite{CBS14}.

In this way, we generate the new nodes in $V_t$ with associated blood types and cPRA rates. In particular, for the latter, we first generate the cPRA interval of a patient and then, their actual cPRA number is drawn from the uniform distribution 
on that interval. Next, we generate the arcs of $A_t$. If donor $d_i$ is blood type compatible with patient $p_j$, then the arc $(i,j)$ is excluded from the exchange graph $G_t$ 
with probability equal to the cPRA rate of the patient $p_j$.\footnote{We use the standard exclusion probability formula
to map cPRA rate to compatability. We remark, however, that cPRA is not the only factor determining compatability and there is no
perfect mapping formula; see Delorme et al.~\cite{DGG22}.}
It follows that, including altruistic donors, there are $84$ types of nodes, distinguished by blood type and cPRA interval.

We remark it is natural to also incorporate a departure rate into the model. Departures may, for example, be due to a donor withdrawing from the program, a patient receiving a kidney from an alternative source, the death of a patient, etc. This is simple to include in a dynamic model, but we have chosen not to do so as we do not have accurate statistics
on which to base the corresponding event probabilities.
Such departures are excluded for every transplant selection algorithm we test and, a priori, there is no reason to suppose that such departures affect one algorithm more than another. 

In each matching period, a transplant selection algorithm is used to select transplant cycles and paths.
The corresponding patient-donor pairs and altruistic donors are then cleared from the market.
Patients that did not receive a transplant remain in the pool. In addition, new patient-donor pairs and altruistic donors
enter the system, according to the prescribed arrival rates, and
a new set of transplant cycles and chains are selected in the  next period. The process then continues and
our interest is in understanding its long-run behaviour.
We want to study how structural properties of the waiting pool vary when different transplant 
selection algorithms are used to select the cycles and chains. 
In particular, we design seven learning algorithms to select weights for the 
exchange graphs. We then study the performance of these algorithms using a large number of full-scale simulations,
and compare their performances against those of the KPD algorithm and the myopic algorithm.

\section{Algorithms and Experiments}\label{sec:algs}

Here, in Section~\ref{sec:learning} we present our learning algorithms and explain how they generate node weights for the exchange graphs. In Section~\ref{sec:KPD} we describe the Canadian Kidney Paired Donation (KPD) program and, in particular, the weighting system used by its transplant selection algorithm. Finally, in Section~\ref{sec:experiments}, we outline our experimental approach for testing the performance of transplant selection algorithms. Specifically, we describe the set of simulations used on algorithms from the three classes of weighting system: learned-weights, expert-defined weights (KPD) and myopic weights.

\subsection{Fair Learning Algorithms}\label{sec:learning}
Our learning algorithms will assign a weight to each node
in the exchange graph. In particular, each patient-donor node is assigned a weight $w_\tau$ depending upon its {\em type} $\tau$.
The type is determined by the blood group of the donor,
the blood group of the patient, and the cPRA group of the patient. Because we consider 5 cPRA groups, there are 80 possible types of patient-donor node.\footnote{As stated, practitioners may also add individual weight adjustments to patient-donor nodes but this is irrelevant to the conclusions of this work. }

Recall that our aim is to obtain a fair algorithm, specifically
equity in terms of transplant accessibility.
To do this, we desire a set of weights {\bf w} that, in the long-run, contains each type in (roughly) the same
proportion in the queue (program) as their proportion in the general population. Ergo, we define ${\tt pop}_\tau$ to be the
proportion of type $\tau$ in the population and
${\tt que}_\tau({\bf w})$ to be the proportion of type $\tau$ in the waiting pool in the long-run using the weighting system {\bf w}.
For our purposes, we define the long-run to be 50 time periods, that is, $16\frac23$ years with matching rounds every four months.

We pre-learn the weights using a sequence of 50 simulations each consisting of the aforementioned 50 matching periods. In each matching period, 
given the weights,
we used the \textit{position-indexed chain-edge formulation} 
by Dickerson et al.~\cite{DMP16} to calculate the optimal
choice of transplant paths and cycles.\footnote{Code, by James Trimble, for the position-index formulations is available at \url{https://github.com/jamestrimble/kidney\_solver}.} 
We let ${\bf w}^t$ be the weights used throughout the 50 periods of the $t$th simulation. To begin we set ${\bf w}= {\bf w}^1$ where $w^1_\tau=1$ for each patient-donor type $\tau$, that is, the myopic weighting system.
At the end of the $t$th simulation we update ${\bf w}$ via
$$w^{t+1}_\tau = f\left(\frac{{\tt que}_\tau({\bf w}^t)}{{\tt pop}_\tau}\right),$$
where $f(x)$ is an update function. Evidently, the update function must have the property that it is monotonically increasing in $x$, that is, in 
$\frac{{\tt que}_\tau({\bf w}^t)}{{\tt pop}_\tau}$.
To accomplish this we define two classes of update function, namely linear updates and exponential updates.
The linear rule, {\tt{Lin}}(a), is defined as 
$$f(x) = 1+ \frac{x}{a}.$$
The exponential rule, {\tt{Exp}}(a), is defined as 
$$f(x) = (a+1) - a e^x.$$
Furthermore, we then scale the functions so that the minimum
weight is at exactly $1$.

For the linear rule we select $a\in \{1, 2\}$ and for
the exponential rule we select $a\in \{1,3,5,7,9\}$.
This produces seven learning algorithms. The performance of each algorithm will be evaluated (via experiments
described below)
using the weights ${\bf w}={\bf w}^{51}$ obtained at the end of the sequence of 
simulations for that algorithm.

In addition to the 80 types of patient-donor node, there are 4 types of altruistic donor node.
Until now the weight of such a node type has been set to $w_\tau=0$. However, we rerun each learning algorithm over a range of different (predominantly, negative) weight values W for the altruistic donors. Similarly, weights were learnt for a range of maximum path and chain lengths and for various altruistic arrival rates.

\subsection{The Canadian KPD Algorithm}\label{sec:KPD}
Here we give a brief overview of the program;
we refer the reader to Cole et al.~\cite{CNC15} and Malik and Cole~\cite{MC14} for detailed descriptions of the program. The KPD program uses an expert-defined weighting system.
The nodes/arcs in the exchange graph are weighted according a point system shown in Table~\ref{tab:kpd}.
We remark that the points (weights) used in the KPD algorithm are of a higher magnitude than those used in our learning algorithms; since this is equivalent to a change in scale, it will have no impact on our subsequent results.

\begin{table}[!h]
\begin{center}
\begin{tabular}{|l|c| } 
\hline
{\bf \qquad \qquad \qquad \qquad Attributes} & {\bf Match Points} \\
\hline\hline
{\em Any Transplant} & 100 \\
{\em Highly Sensitized (cPRA $\ge 0.80$)} & 125 \\
{\em ABO Match: O to O} & 75 \\
{\em Paediatric Recipient} & 75 \\
{\em Recipient is Prior Living Donor} & 75 \\
{\em HLA 0/6 Mismatch} & 75 \\
{\em Dialysis Wait Time (starting at initiation on dialysis)} & Days $\div$ 30 \\
{\em Geography: Same City} & 25 \\
{\em Donor/Recipient Age $\Delta\le 30$ years} & 5 \\
{\em ABO Match: A to A, B to B, AB to AB} & 5 \\
{\em EBV Negative to Negative Match} & 5 \\
\hline
\end{tabular}
\caption{The KPD point system~\cite{MC14}.}\label{tab:kpd}
\end{center}
\end{table}

\begin{figure}[h!]
\begin{center}
\includegraphics[scale=0.65]{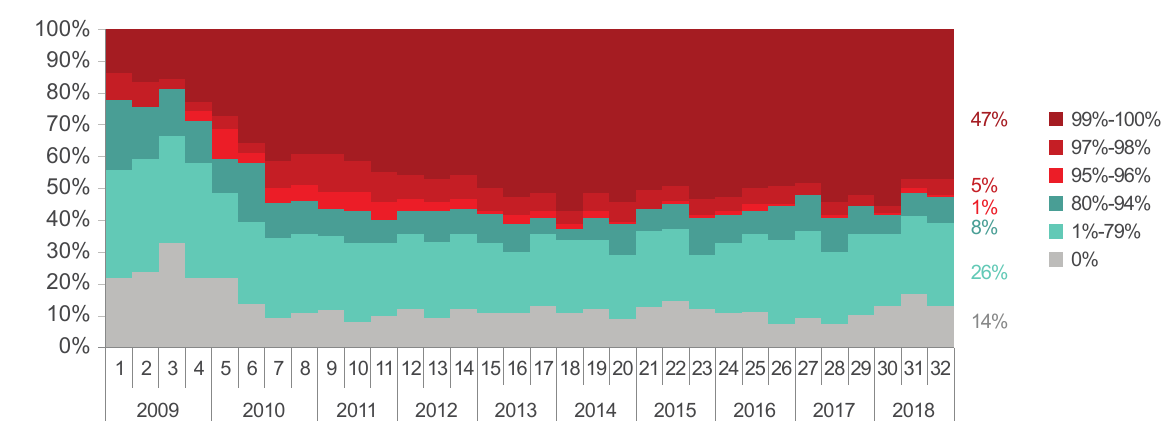}
\end{center}
\caption{KPD pairs by cPRA of patients over times~\cite{CBS18}.}\label{fig:cpra}
\end{figure}

Interestingly, this is a revised point system that was devised in part to address the accumulation of hard-to-match patients in the KPD waiting pool. This is illustrated in Figure~\ref{fig:cpra}. 

Indeed, the program has been quite successful in stabilizing the proportion of hard-to-match patients in the pool.
Care should be taken with Figure~\ref{fig:cpra}, however,
as it indicates only the proportions and {\bf not} the magnitudes of patients types.
This observation is key in that limiting the growth of the waiting pool is fundamental in the performance of a transplant selection algorithm.

Transplant selection occurs every 4 months in the KPD program. It is important to note that the matches proposed by the KPD algorithm (or, indeed, by any transplant selection algorithm) may not all be undertaken, e.g., due to surgical concerns or the withdrawal of a patient. 
To date, such factors cannot be accurately modelled. Furthermore, there are some aspects of the KPD scoring scheme that, whilst very simple to implement, we have excluded because we do not have the necessary statistical data to accurately incorporate them into our learning models. These include
donor/recipient ages, geographical location, and time on dialysis. Regardless, a priori, there is no reason to suppose that these excluded elements will affect one transplant selection algorithm more than another; thus, for the purposes of algorithmic comparison we assume
it is safe to omit these factors from our model.

Finally, we remark that for the KPD program the surplus donor at end of a path is transferred to the standard kidney donation program for unpaired patients as a donor for the medical centre that registered the altruistic donor~\cite{CNC15}, that is, as a living donor on the deceased donor list.

\subsection{The Experiments}\label{sec:experiments}
So we have described how weights are calculated for the 
mypoic algorithm, the KPD algorithms, and our learning algorithms. Given these weighting systems how do the respective
algorithms perform?
We tested this using a large number of experiments run via the dynamic kidney exchange model.
That model, described in Section~\ref{sec:model}, was implemented using Python 3. As stated, for these experiments, simulated patients and donors were generated using Canadian population data statistics based on publicly available data from Canadian Blood Services. 
For each algorithm, we conducted 50 random simulations, each lasting for 50 matching rounds, that is, over a time-span of 200 months.
In each matching period, 
given the weights prescribed by the relevant algorithm,
we again use the \textit{position-indexed chain-edge formulation} 
by Dickerson et al.~\cite{DMP16} to calculate the optimal
choice of transplant paths and cycles. These experiments were carried out over a range of maximum cycle lengths ($C\le 5$) and maximum path lengths ($P\le 20$).
Using these experiments we evaluate the long-term impact of each of these transplant selection algorithms.
Since our algorithms are designed for group fairness, before describing the results of our learning experiments in Section~\ref{sec:results} let us first discuss how to measure fairness in kidney exchange programs.


\section{Fairness in Kidney Exchanges}\label{sec:fair}

We will present the results of our simulations in Section~\ref{sec:results}. These compare the performances of a variety of 
kidney exchange algorithms (i.e., transplant selection objectives) with respect to the number of transplants, waiting times, the utilization rate of altruistic donors, and fairness.
The meanings of the first three of these criteria are self-evident. But what do we want from a ``fair'' kidney-exchange algorithm?
We address this question in this section and, as a result, provide a class of measures of fairness by which we will compare
our algorithms.

\subsection{Measures of Fairness}\label{sec:measures}
An extremely important desideratum for kidney exchange mechanisms is fairness across groups and individuals.
The first measures of fairness in kidney exchanges were provided by Dickerson et al.~\cite{DPS14}.
Motivated by the {\em price of fairness} for allocation problems~\cite{BFT11,CKK12} their
fairness measures (lexiographic and weighted) concerned fairness across groups defined by cPRA rates.
Fairness in kidney exchanges is now well-studied in the literature, typically based on measures of group 
fairness. A few measures based on individual fairness have been introduced; see Farnadi et al.~\cite{FAB21} and St-Arnaud et al.~\cite{St22}. Furthermore, the provision of fairer allocations via the incorporation of human values into matching algorithms has
also been proposed; see, for example, Freedman et al.~\cite{FBS20}.

There is a very natural way to evaluate group fairness using the weighted 
power mean. Specifically, given $k$ groups with weights $\pmb{\alpha}=(\alpha_1, \alpha_2, \dots, \alpha_k)$ and utilities
${\bf u}=(u_1, u_2, \dots, u_k)$, we may define group utility by the {\em weighted power mean}
\begin{equation}\label{eq:wpm}
M^{\rho}_{{\bf \alpha}}({\bf u}) \ =  \left( \sum_{i=1}^k \alpha_i\cdot u_i^\rho \right)^{\frac{1}{\rho}}.
\end{equation}
We refer the reader to the book~\cite{Bul03} by Bullen for technical details on weighted power means and their mathematical properties. Important here
is the fact that~(\ref{eq:wpm})
corresponds to the well-studied {\em constant elasticity of substitution} (CES) utility functions in economics~\cite{Sol56},\cite{ACM61},\cite{DS77}, where the
elasticity of substitution is $\sigma=\frac{1}{1-\rho}$. What is interesting about this class of functions is that $\rho$ can be adjusted to 
give a variety of fairness measures over the groups. Most important are the special cases $\rho \in \{1, 0, -\infty\}$.
The case $\rho=1$ corresponds to {\em weighted utilitarian social welfare}, where we wish to maximize $\sum_{i=1}^k \alpha_i\cdot u_i$.
The case $\rho=0$ corresponds to {\em weighted Nash social welfare}, where we wish to maximize $\prod_{i=1}^k u_i^{\alpha_i}$, also known as
the Cobb-Douglas utility. Finally, the case $\rho=-\infty$ corresponds to {\em egalitarian (Rawlsian) social welfare}, where we wish to 
maximize $\min_{i=1}^k u_i$, also known as the Leontif utility. 

Thus (\ref{eq:wpm}) induces a class of fairness measures. However, it still remains to define the groups, the weights and utilities. Let us begin with the groups. Recall a major motivation for this work is the accumulation of hard-to-match patients in the waiting pool,
Now the cPRA score estimates the proportion of 
donors against which the patient has antibodies rendering the match incompatible; thus, the higher the cPRA rate the
harder it is to match a patient. Consequently, as in~\cite{DPS14}, it is natural to group the patients according to cPRA rate.
We will use five groups for cPRA rates in the intervals 
$(0, 0), (0.01, 0.50), (0.51, 0.94), (0.95, 0.96)$ and $(0.97, 1)$. 
We may then set the weight of a group to be the probability a patient lies in that group, that is 
$(\alpha_1, \alpha_2, \alpha_3, \alpha_4, \alpha_5)=(0.24, 0.29, 0.24, 0.10, 0.13)$. 
Finally, we need to define the utility $u_j$ of a group $j$. This utility should be inversely related to the accumulation of the group in the patient pool, where
we denote by $q_j$ the number of patients of group $j$ in the pool at the end of the experiments. 
In addition, we desire that the utility scales linearly; specifically if the quantity $q_i$ of every group
increases by a factor $c$ then the utility of each group should decrease by a factor~$c$. 
These properties can be obtained by defining the utility of a group to be its weight divided by the number of its members in the queue.
That is, $u_j =\frac{\alpha_j}{q_j}$. We can view this utility in the following way. Let $Q=\sum_{i=1}^k q_i$ be the total queue size, 
and let $\beta_j=\frac{q_j}{Q}$ be the proportion of the queue made up of members of group $j$. 
Then $u_j= \frac{1}{\beta_j/\alpha_j}\cdot \frac{1}{Q}$. Thus the utility of the group decreases as the group becomes over 
represented in the queue (when $\beta_j$ grows relative to $\alpha_j$) and as the queue size $Q$ increases. 

So, with these weights and utilities, the weighted power mean (\ref{eq:wpm}) induces three fairness scores for $\rho \in \{1, 0, -\infty\}$ 
by which we will evaluate our kidney exchange algorithms in Section~\ref{sec:results}.
These scores illustrate varying trade-offs between the utilities of each group.
At one extreme, utilitarian social welfare $(\rho=1)$, the objective is maximize collective welfare; all groups are important
but specific groups may be penalized if this leads to greater benefits for the other groups.
At the other extreme, egalitarian social welfare $(\rho=-\infty)$, the objective is to maximize the utility of the worst-off group.
For the purpose of KEPs, egalitarian social welfare is the most appropriate fairness objective with respect to the 
accumulation of hard-to-match patients in the waiting pool. So this measure will be our primary focus.
However, Nash social welfare $(\rho=0)$ is also an extremely useful fairness measure for allocation problems (such as KEPs).
This is highlighted by Caragiannis et al.~\cite{CKM19} in their influential work on fairness. They state ``the Nash social welfare solution exhibits an elusive 
combination of fairness and efficiency properties, and can easily be computed in practice. It provides the most practicable approach to 
date -- arguably, the ultimate solution, for the division of indivisible goods under additive valuations''.
Informally, located between utilitarian and egalitarian social welfare, Nash social welfare encourages improvements in collective welfare whilst, 
at the same time, penalizing imbalances between groups.

%
%
%
%

\section{Results}\label{sec:results}

We now present the results of our experiments.
In Section~\ref{sec:summary} we provide a summary of
the performance of the transplant selection algorithms for a variety of criteria including the number of transplants, waiting times and group fairness.
These results highlight the value of learned-weights. 
Then we provide more in-depth analyses. 
First, in Section~\ref{sec:paths}, we explore the sensitivity of the transplant selection algorithms to cycle and path length caps.
Then, in Section~\ref{sec:altruists}, we show the vital importance of assigning negative weights to altruistic donors. Finally, in Section~\ref{sec:metrics}, we detail results on performance with respect to our three group fairness metrics.

\subsection{Summary of Results}\label{sec:summary}
Table~\ref{tab:main} shows the performance of the KPD and myopic algorithms, and the seven learning algorithms. There are four types of criteria for comparison.
The three main criteria are:
\begin{enumerate}
\item {\em Number of Transplants}. This is shown in two equivalent ways. First the
total number of transplants\footnote{The surplus donor at the end of a path may subsequently donate on the deceased donor program. Such matches are not counted in the total number of transplants.} achieved over the 50 time periods by each algorithm, on average over all experiments.
Second, the corresponding percentage of patients who received a transplant.
\item {\em Waiting Times.}  This is measured in two different ways. The first is the average wait time of each transplant recipient in months.
This is an important measure but it does not give the entire picture as it excludes those patients that remain on the waiting list, in particular, harder to match patients.
Thus, the second way is the average wait time, which includes the waiting times of every patient.
\item {\em Group Fairness}. We measure fairness in three ways: utilitarian group fairness, Nash group fairness and
egalitarian group fairness. Only our primary fairness measure, egalitarian group fairness, is shown in Table~\ref{tab:main}. The measure is scaled so that the KPD 
algorithm has measure $1$ and the higher the measure the fairer the algorithm.
(Detail discussion and our results concerning all three fairness
measures will follow.)
\end{enumerate} 
The significance of the fourth criteria, {\em ``Altruist Usage''}, (and that of the negative weight, $W$, assigned to each altruistic donor as in Table~\ref{tab:main})
will be shown in Tables~\ref{tab:altruist-weights} and~\ref{tab:prices}.

\begin{table}[h!]
\begin{center}\fontsize{9}{11}\selectfont
\setlength\tabcolsep{2pt}
\begin{tabular}{|c||c|c|c|c|c|c|c|} 
\hline
\multicolumn{7}{|c|}{Periods: 50 \ \ \# Participants: 1994.78 \ \  $\lambda_A$: 4.5625  \quad C: 5 \ \  P: 5}\\
\hline\hline
\multirow{2}{2em}{Alg.} 
& \multirow{2}{*}{W}
& \multirow{2}{*}{\#Matches} 
& \multirow{2}{*}{\%Match} 
& \multirow{2}{4.8em}{Wait Time\\(Recipients) } 
& \multirow{2}{2.2em}{Wait\\ Time} 
& \multirow{2}{3.4em}{Egal.\\Fairness}\\
&&&&&&\\
\hline\hline
{\tt Myop} 			&0.0			&1613.76	&80.90	&13.55 &29.52 & 0.60\\ \hline
{\tt KPD}				&0.0			&1612.64	&80.84 	&13.58 &29.56& 1.00\\ \hline
{\tt Lin}(1)				&-2.0			&1698.22	&85.13	&11.51 &23.45& 1.06\\ \hline
{\tt Lin}(2)				&-3.0			&1700.38	&85.24	&11.47 &23.45& 0.85\\ \hline
{\tt Exp}(1)			&-1.5		&1707.46	&85.60	&11.88 &22.62& 1.10\\ \hline
{\tt Exp}(3)			&-2.0			&1694.72	&84.96	&11.69 &23.56& 1.17\\ \hline
{\tt Exp}(5)			&-3.0			&1697.62	&85.10	&11.51 &23.55& 0.91\\ \hline
{\tt Exp}(7)			&-4.0			&1705.88	&85.52	&10.89 &22.86& 0.61\\ \hline
{\tt Exp}(9)			&-5.0			&1697.52 	&85.10	&11.45 &23.61& 1.04\\ \hline
\end{tabular}
\caption{Performances of the transplant algorithms.}\label{tab:main}
\end{center}
\end{table}



The main observation from Table~\ref{tab:main} is that all seven of the learning algorithms outperform the myopic algorithm for every single 
criteria (including utilitarian and Nash group fairness, not shown in the table). 
Four of the learning algorithms outperform the KPD algorithm for every single 
criteria; the remaining three learning algorithms outperform the KPD algorithm on every criteria
except egalitarian fairness.
The performance of each learning algorithm is comparable but, overall, {\tt Exp}(1) and {\tt Exp}(3) perform the best. 
The headline statistics are the following. In comparison to the KPD and myopic algorithms, the learning 
algorithm {\tt Exp}(1) leads to a 6\% increase in the number of transplants (from about 1613 to 1707), 
a decrease of 24\% in the average waiting time (from about 29.5 to 22.6 months) and an improvement of 10\% and 21\%, respectively, in the egalitarian fairness measure. Furthermore, both {\tt Exp}(1) and {\tt Exp}(3) reduce by 2\% the use of altruistic donors (see subsequent discussion).


Before presenting more detailed results and analyses, it is worth discussing whether even better algorithms
may be obtainable. In fact, it will be hard to obtain better algorithms. To see this, 
consider the number of transplants. The \%Match column shows the percentage of patients that receive transplants.
The learning algorithm {\tt Exp}(1) achieves 85.6\%. This does not sound spectacular, but it is near optimal
because, when we take these measurements after 50 periods, it is not realistic to provide transplants for 
all the recent arrivals. Indeed, it appears that about a 90\% match rate is the maximum achievable.
(An indication of this fact will be seen in Table~\ref{tab:altruist-number} where a 90\% match rate is achieved 
but only after tripling the number of altruistic donors in the KEP.)


\subsection{Cycle and Path Lengths}\label{sec:paths} The results shown in Table~\ref{tab:main} are based on allowing cycles of length at most $C=5$ and paths of
length at most $P=5$. These restrictions are justified by the following two lessons.\\[-0.3cm]

\noindent {\tt Lesson I.} {\em There is negligible hypothetical benefit in allowing for paths of length greater than 5.}\\[-0.3cm]

\noindent {\tt Lesson II.} {\em It is very important to allow for cycles of length 3. The hypothetical improvements 
with cycles length of 4 or 5 are marginal.}\\[-0.3cm]

We remark that these two lessons will not surprise practitioners. 
Whilst the potential importance of long paths was originally highlighted by~\cite{AGRR12}, 
early works gave theoretical and experimental evidence that consideration of only short cycles~\cite{RSU07} or short paths~\cite{DPS12b} suffice.

We begin with the justification of the cycles cap. 
In kidney exchange algorithms, cycle lengths are restricted for two major reasons. First, it is generally desirable
that transplant operations for cycles are carried out simultaneously. Specifically, it ends the risk of a donor
dropping out after their accompanying patient receives a transplant -- this renders the cycle incomplete,
leaving another patient without a transplant despite their accompanying donor having donated a kidney.
For logistical reasons, e.g. the number of transplant teams, it is then impossible to deal with long cycles. 
Second, computing long cycles is computationally more difficult. 
Given the huge number of integer programming problems (one for each application of a kidney-exchange algorithm) solved in computing weights by our learning algorithms and the large number of tests 
required by our experiments, we found that $C=5$ was the best choice.\footnote{For example, the case
of $C=6$ can be solved in a reasonable amount of time and is thus feasible for use in a kidney exchange program.
However, it is not currently reasonable for repeated usage to the extent required for our tests.
Regardless, our results suggest solving for $C=6$ would provide little benefit.}

On the other hand, these two problems do not arise with paths. First, there is no need for the operations to take place simultaneously; a donor may still drop out, but since the path began with an altruistic donor there is no 
accompanying patient left in the lurch.\footnote{Of course, if transplants along a path are conducted sequentially, not simultaneously, then the time span between matching rounds will impose a natural restriction on the maximum path length.}  Second, computationally searching for solutions with long paths is not 
prohibitive even with the mass repetition involved here. But ``in practice, it is desirable
to impose a cap on chain length, since there is an increasing chance that the final
exchanges planned for a chain will not proceed as the chain length is increased (due to
various reasons such as pre-transplant crossmatch incompatibility, death of a recipient
before transplant, the recipient receiving a deceased-donor kidney, and so on)''~\cite{DMP16}.
Fortunately, our simulations show that there is no hypothetical benefit in allowing for long paths.
\begin{table}[h!]
\begin{center}\fontsize{9}{11}\selectfont
\begin{tabular}{|c|c|c|c|c|c|c|c|c|} 
\hline
\multicolumn{9}{|c|}{{\tt Exp(3)} \quad Periods: 50 \quad \# Participants: 1994.78 \quad $\lambda_A$: 4.5625  \quad C: 5 \quad W:-2}\\
\hline\hline
\multirow{2}{1em}{P}& 
\multirow{2}{4.5em}{\#Matches} &
\multirow{2}{4em}{\%Match} & 
\multirow{2}{4em}{Altruist \\\% Usage}&
\multicolumn{4}{|c|}{\# Paths}&
\multirow{2}{4.5em}{\% of Path Matches} \\
&&&&1-5&6-10&11-15&16-20&\\\hline
5	&1698.06	&85.13	&97.25 	&228.28	&0.00		&0.00		&0.00		&50.24 \\ \hline
10	&1704.20	&85.43	&96.99 	&114.18	&83.46	&0.00		&0.00		&63.91\\ \hline
15	&1704.20	&85.43	&97.07  	&159.74	&43.76	&21.36	&0.00		&61.65\\ \hline
20	&1701.16	&85.28	&97.10  	&166.06	&37.12	&12.46	&7.19	&59.60\\ \hline
\end{tabular}
\caption{The impact of constraining path lengths.}\label{tab:paths}
\end{center}
\end{table}

Specifically, Table~\ref{tab:paths} shows how the expected number of matches over 50 periods changes
when we allow paths to have lengths up to $5, 10, 15$ or $20$. We immediately see that the increase in
matches is negligible (from 85.13\% to 85.43\%) when we allow for paths of length greater than $5$. 
We remark that Table~\ref{tab:paths} concerns the learning model {\tt  Exp(3)}. Here, we chose to illustrate our results with this model due to its high performance and because
it was the best for the egalitarian fairness criterion, but a similar picture arises with all the models.
Consequently, for our experiments we imposed a maximum path length of $P=5$, and we conclude {\tt Lesson I}.




\begin{table}[h!]
\begin{center}\fontsize{9}{11}\selectfont
\begin{tabular}{|c|c|c|c|c|c|c|c|c|c|} 
\hline
\multicolumn{10}{|c|}{{\tt Exp(3)} \quad Periods: 50 \quad \# Participants: 1994.78 \quad $\lambda_A$: 4.5625  \quad P: 5 \quad W:-2}\\
\hline\hline
\multirow{2}{1em}{C}& \multirow{2}{4.5em}{\#Matches} &\multirow{2}{4em}{\%Match} &\multicolumn{4}{|c|}{\# Cycle}
&\multirow{2}{3em}{\# Paths} &\multirow{2}{4.8em}{\# Patients\\ in Paths} &\multirow{2}{4.5em}{\% of Path Matches} \\
&&&2&3&4&5&&&\\\hline
2	&1631.22	&81.77	 	&263.40	&0.00		&0.00		&0.00   		&229.28   &853.14		&67.70 \\ \hline
3	&1685.34	&84.49	 	&95.16	&162.60	&0.00		&0.00 		&228.34   &894.64 		&59.76\\ \hline
4	&1690.40	&84.74		&58.44	&70.88	&116.56	&0.00		&219.86	&1007.22		&52.92\\ \hline
5	&1698.06	&85.13	 	&48.14	&46.12	&53.12	&79.56 	&222.28 	&1104.42		&50.24\\ \hline
\end{tabular}
\caption{The impact of constraining cycle lengths.}\label{tab:cycles}
\end{center}
\end{table}
Next, does it actually help to have cycle lengths up to $5$? This question is addressed in Table~\ref{tab:cycles}.
Observe that there is a very significant improvement when we increase the maximum cycle length 
from $2$ to $3$. Specifically, the percentage of matches increases from 81.8\% to 84.5\%.
However, increasing the cycle length beyond that leads only to a small improvement, up to 85.13\%
when cycles of length $5$ are allowed;
this small hypothetical improvement in the number of transplants may be counteracted by the practical disadvantages inherent in using longer cycles. Hence, we obtain {\tt Lesson II}.

\subsection{Altruist Weights}\label{sec:altruists}
Our learning models learn positive weights for patient-donor nodes in the exchange graph. 
However, it turns out to be vital to also assign a negative node weight to an altruistic donor.
Before investigating the effect of this modification, let us understand why it may be useful.
The incorporation of positive weights on patient-donor nodes signifies the ``value''
of performing those transplants. In contrast, a negative weight can be interpreted as signifying the
``price'' of using an altruistic donor. An altruistic donor can be viewed as a scarce and valuable resource 
in an kidney exchange mechanism. They must be used prudently and the standard economic approach
to enforce this is to assign them a price -- the consequence of this price means an altruistic donor will only be used 
when the resultant benefit is significant.
The reader may ask whether this is really helpful. The answer is {\tt yes}: there is a dramatic difference in the 
performance of the kidney-exchange algorithm (in terms of the number of transplants, waiting times, and fairness) with
different altruistic donor weights. 
Table~\ref{tab:altruist-weights} shows this emphatically for the learning model {\tt Exp}(3). 
When the altruist weight is zero (or even positive) the number
of matches is low, around 80.7\%. As the weight falls, the number of matches increases dramatically
to around 85\% for the weight range $[-2,-3]$. Further decreases in weight beyond $-3$ then lead to a large 
fall in the number of transplants.
\begin{table}[h!]
\begin{center}\fontsize{9}{11}\selectfont
\begin{tabular}{|c|c|c|} 
\hline
\multicolumn{3}{|c|}{{\tt Exp(3)}  Periods: 50 \ \ \# Participants: 1994.78 \  C: 5   P: 5}\\
\hline\hline
W & Altruist \% Usage  &   \%Match \\
\hline\hline
0.00			&99.64	&80.77	\\ \hline
-0.75			&98.21	&83.69	\\ \hline
-1.50			&97.71	&84.57	\\ \hline
-2.00			&97.32	&84.96	\\ \hline
-2.50			&97.34	&85.42	\\ \hline
-3.00			&96.12	&85.18	\\ \hline
-5.00			&51.90	&81.01	\\ \hline
-10.00			&40.23	&75.45	\\ \hline
-15.00			&0.00	&73.83	\\ \hline
\end{tabular}
\caption{The effect of varying the weight of altruistic donors.}\label{tab:altruist-weights}
\end{center}
\end{table}


This is very natural. When the price is under $2$ it is too low and the demand for altruists is too high. The result is 
that the altruists are used in a suboptimal manner, for example, to initiate short paths with easy-to-match patients.
When the price is above $3$ it is too high and the demand for altruists is too low. The   
result is that the potential benefits of using the altruists are not all achieved.
This then explains why a negative weight $W$ for the altruists is given in Table~\ref{tab:main} for each model. 
We remark that the optimal price for an altruistic varies relative to the other weights given by the model for the arcs 
or patient-donor pairs. Hence, it also explains why in Table~\ref{tab:main} the best choice of altruist weight
differs for each model.
In this way, we obtain our most important lesson: \\[-0.3cm]

\noindent {\tt Lesson III.} {\em Assigning an appropriate negative weight (price) to each altruistic donor is fundamental 
in the determining the quality of outcomes.}\\[-0.3cm]

To see why this lesson is critical, let us see what happens when we apply it to the myopic and KPD algorithms. 
Concretely, assigning a price to the altruists gives the results shown in Table~\ref{tab:prices}, where we call 
the modified algorithms Myopic$^+$ and KPD$^+$.
This trivial modification has a massive impact. Both Myopic$^+$ and KPD$^+$ now
perform a comparable number of transplants to our learning algorithms 
roughly 85\% each (up from 81\%). There is also a very large improvement in fairness. The egalitarian group fairness 
measure increasing by 26\% for the KPD algorithm and by 77\% for the myopic algorithm.
Finally, waiting times also improve significantly with one exception: for the KPD algorithm, the average waiting time decreases
from 29.6 to 23.6 months but waiting times per recipient increase from 13.6 to 14.8 months. 
The explanation for this counter-intuitive behaviour is that KPD$^+$ performs many more transplants
than KPD but many of these additional operations involve hard-to-match patients who have naturally spent
longer than average time on the waiting lists.  

\begin{table}[h!]
\begin{center}\fontsize{9}{11}\selectfont
\setlength\tabcolsep{5pt}
\begin{tabular}{|c||c|c|c|c|c| } 
\hline
\multicolumn{6}{|c|}{Periods: 50 \ \ \# Participants: 1994.78 \ $\lambda_A$: 4.5625 \ C: 5  \ P: 5}\\
\hline\hline
\multirow{2}{*}{Alg.}
& \multirow{2}{*}{W}
& \multirow{2}{*}{\%Match}
& \multirow{2}{5em}{Wait Time\\ (Recipients)} 
& \multirow{2}{2em}{Wait\\ Time} 
& \multirow{2}{2.8em}{Egal.\\ Fairness}\\
&&&&&\\
\hline\hline
{\tt Myop} 			&0.0			&80.90	&13.55 &29.52 & 0.60\\ \hline
{\tt Myop+} 			&-2.0			&84.80	&12.86 &23.87 & 1.06\\ \hline
{\tt KPD}				&0.0		&80.84 	&13.58 &29.56 & 1.00\\ \hline
{\tt KPD+}				&-150.0		&85.07 	&14.82 &23.64 & 1.26\\ \hline
{\tt Exp}(1)			&-1.5		&85.60	&11.88 &22.62 & 1.10\\ \hline
{\tt Exp}(3)			&-2.0			&84.96	&11.69 &23.56 & 1.17\\ \hline
\end{tabular}
\caption{The impact of pricing altruistic donors.}\label{tab:prices}
\end{center}
\end{table}

Consequently, excluding waiting time per recipient, the KPD$^+$ algorithm has performance 
comparable to our learning algorithms. An important conclusion can be drawn from this: the expert-designed 
weights in the KPD system work very well; the inferior performance of the KPD system is almost entirely due to the lack of 
an altruistic donor price!
A second important conclusion can be drawn by comparing Myopic$^+$ with KPD. We have seen
that the choice of KPD weights are appropriate. However, Myopic$^+$ significantly outperforms
KPD in all three measures. But Myopic$^+$ is an unweighted system except for the 
altruistic donor price. Consequently, choosing a (negative) weight for the altruists
is more effective than assigning an accurate (positive) weight for {\em every} other patient-donor pair. 

Given the huge impact of altruistic donors, a natural question is what would be the performance
of a kidney exchange if the number of altruistic donors changed.
We examine this in Table~\ref{tab:altruist-number} where the altruist arrival rate $\lambda_A$ varies from $0$ to $10.5$.
Recall, in the Canadian KPD program the current altruist arrival rate is roughly $\lambda_A=4.5$.
\begin{table}[h!]
\begin{center}\fontsize{9}{11}\selectfont
\begin{tabular}{|c|c|c|c|c|c|} 
\hline
\multicolumn{6}{|c|}{{\tt Exp(3)} Periods: 50 \ \# Participants: 1994.78 \ \   C: 5   P: 5  W:-2}\\
\hline\hline
\multirow{2}{2em}{$\lambda_A$}
&   \multirow{2}{4em}{\%Match} 
& \multirow{2}{4em}{Matches/ \\Altruist} 
&\multirow{2}{*}{\#Altruists}
& \multirow{2}{4em}{\#Paths}
& \multirow{2}{4em}{Donors/\\ Altruist}\\
&&&&&\\
\hline\hline
0.0			&73.74	&- 	&0.00		&0.00		&-  \\ \hline
1.5			&78.02	&1.14 	&76.88	&75.24	&0.98\\ \hline
3.0				&81.47	&0.92 	&150.12	&146.38	&0.98\\ \hline
4.5			&85.03	&0.95 	&224.92	&219.94	&0.98\\ \hline
6.0			&87.85	&0.75 	&299.16	&281.62	&0.94\\ \hline
7.5			&89.14	&0.34 	&371.48	&311.86	&0.84\\ \hline
9.0			&89.71	&0.15	&448.96	&324.96	&0.72\\ \hline
10.5			&90.04	&0.09 	&521.76	&335.50	&0.64\\ \hline
\end{tabular}
\caption{The impact of the altruistic donors arrival rate.}\label{tab:altruist-number}
\end{center}
\end{table}

Table~\ref{tab:altruist-number} is extremely informative. For example it shows that, 
for the current altruist arrival rate $\lambda_A=4.5$, each altruistic donor saves almost two lives. To see this note that
each new altruist increases the number of transplants by $0.95$ on average. In addition, the 
altruist usage rate is 98\%. But the use of an altruist at the start of a path means the donor at the end of the path is 
unused. In the Canadian KPD, such donors are currently requested to then donate in the deceased donor program.
This means each altruistic donor then creates on average an additional $0.98$ transplants in deceased donor program.
Thus $0.95+0.98=1.93$ is a good measure of the number of lives saved by each altruistic donor. 

Instead, such a donor could be requested to remain in the paired donation program
as a {\em defacto} altruistic donor, called a ``bridge donor''. 
Observe that the use of a bridge donor recursively then generates another bridge donor.
Consequently, an important open question is to determine whether this {\em multiplier effect}
implies the benefits of the donor remaining in the paired donation program are greater than
the benefits of the donor moving to the deceased donor program.

Table~\ref{tab:altruist-number} also shows that the marginal increase in transplants per altruist is {\em decreasing}. For example, 
between $\lambda_A=0$ and $\lambda_A=1.5$ each new altruist increases the number of transplants by $1.14$ on average. 
The marginal increase remains around $1$ until $\lambda_A=6$ and then drops precipitously from $\lambda_A=7.5$ upwards.
Furthermore, when the altruist arrival rate becomes very high, the altruist usage rate drops noticeably.
Indeed for $\lambda_A\ge 9$ the number of lives saved per altruist falls below $1$.
For example, for $\lambda_A=10.5$ marginal increase in transplants is just $0.09$ and the altruist usage rate is just 64\%. 
This lead us to our next lesson: \\[-0.3cm]

\noindent {\tt Lesson IV.} {\em The marginal contribution in terms of lives saved per altruistic donor is non-increasing
in the number of altruists. (At the current altruist arrival rate in Canada, each altruistic donor saves approximately 
two lives.)} \\[-0.3cm]

Of course this lesson extends beyond the Canadian KPD program. For any kidney exchange program, it is important
to determine how much more effective the program would be with additional altruistic donors.
In particular, the key factor is where on the marginal contribution curve the current altruist arrival rate lies.
We set aside the question related to the ethic of encouraging altruistic donors; see Patel et al.~\cite{PCP11} and Woodle et al.~\cite{WDA10} for discussions.

\subsection{Fairness}\label{sec:metrics}
In Section~\ref{sec:fair} we defined a class of welfare functions to evaluate group fairness using 
the weighted power mean. This class includes three important cases: utilitarian social welfare, Nash social welfare and
egalitarian social welfare.
To compute these welfare functions it is first necessary, for each of our transplant selection algorithms, to 
evaluate the waiting list distribution at the end of our experiments with respect to each group. 
These results are shown in Table~\ref{tab:queues}; we use the altruistic donor weights given in Table~\ref{tab:main}.
\begin{table}[h!]
\begin{center}\fontsize{9}{11}\selectfont
\begin{tabular}{|c||c|c|c|c|c|} 
\hline
\multicolumn{6}{|c|}{Periods: 50 \ \ \# Participants: 1994.78  \quad C: 5 \quad  P: 5}\\
\hline\hline
&\multicolumn{5}{|c|}{Queue Length per cPRA Group}\\
\hline
Model &G1&G2&G3&G4&G5\\
\hline\hline
{\tt Myopic} 			&47.78	&72.14	&89.24	&55.42	&118.42	\\ \hline
{\tt Myopic+} 			&30.62	&48.10	&65.14	&47.02	&111.18	\\ \hline
{\tt KPD}				&88.22	&121.28	&87.66	&14.50	&70.48	\\ \hline
{\tt KPD+}				&66.90	&94.36	&71.90	&10.86	&57.78	\\ \hline
{\tt Lin}(1)				&105.00	&31.22	&75.96	&17.66	&66.72	\\ \hline
{\tt Lin}(2)				&24.24	&61.15	&111.48	&14.16	&82.96	\\ \hline
{\tt Exp}(1)			&19.36	&129.38	&33.56	&41.16	&63.86	\\ \hline
{\tt Exp}(3)			&30.02	&88.52	&111.31	&13.58	&56.62	\\ \hline
{\tt Exp}(5)			&38.54	&114.56	&25.16	&59.34	&59.56	\\ \hline
{\tt Exp}(7)			&6.84	&25.76	&78.74	&62.80	&114.76	\\ \hline
{\tt Exp}(9)			&20.56	&126.42	&37.08	&45.50	&67.70	\\ \hline
\end{tabular}
\caption{Queue Lengths.}\label{tab:queues}
\end{center}
\end{table}

 This produces the social welfare scores shown in Table~\ref{tab:fairness}.
In addition, for the purposes of easy comparison, we create three group fairness measures by scaling the
welfare scores by a fixed constant such that the score of the KPD algorithm is exactly one.
This induces three group fairness measures: utilitarian, Nash and
egalitarian group fairness. 

\begin{table*}[h!]
\begin{center}\fontsize{9}{11}\selectfont
\begin{tabular}{|c||c|c|c||c|c|c| } 
\hline
\multicolumn{7}{|c|}{Periods: 50 \quad Avg Total Participants: 1994.78 \quad $\lambda_A$: 4.5625  \quad C: 5 \quad  P: 5}\\
\hline\hline
&\multicolumn{3}{|c||}{Group Welfare Score}&\multicolumn{3}{|c|}{Group Fairness Measure}\\
\hline
Model &Utilitarian&Nash&Egalitarian&Utilitarian& Nash& Egalitarian\\
\hline\hline
{\tt Myopic} 				&0.00334	&0.00300	&0.00110	&1.14	&  1.10	&0.60 \\ \hline
{\tt Myopic+} 				&0.00488	&0.00415	&0.00117	&1.66	&  1.52	&0.64\\ \hline
{\tt KPD}					&0.00293	&0.00274	&0.00184	&1.00		&  1.00		&1.00     \\ \hline
{\tt KPD+}					&0.00377	&0.00349	&0.00225	&1.28	&  1.27 	&1.22\\ \hline
{\tt Lin}(1)					&0.00482	&0.00398	&0.00195	&1.64	&  1.45	&1.06\\ \hline
{\tt Lin}(2)					&0.00517	&0.00421	&0.00157	&1.76	&  1.54 	&0.85\\ \hline
{\tt Exp}(1)				&0.00585	&0.00444	&0.00204	&1.99	&  1.62	&1.10\\ \hline
{\tt Exp}(3)				&0.00442	&0.00380	&0.00216	&1.51	&  1.39	&1.17\\ \hline
{\tt Exp}(5)				&0.00497	&0.00407	&0.00169	&1.69	&  1.49	&0.91\\ \hline
{\tt Exp}(7)				&0.01272	&0.00659	&0.00113	&4.34	&  2.41	&0.61\\ \hline
{\tt Exp}(9)				&0.00549	&0.00423	&0.00192	&1.87	&  1.55 	&1.04\\ \hline
\end{tabular}
\caption{Measures of Fairness.}\label{tab:fairness}
\end{center}
\end{table*}


To understand this table, recall that all three measures have the property that if the number of patients in {\em every}
group change by an identical factor then the group fairness measure will change by the same factor. For example,
if the number of patients in the waiting list falls by 10\% for {\em every} group then each measure will increase by 10\%.
However, the three measures differ in how they penalize imbalances, with extreme penalties
in the case of egalitarian group fairness.

For the both utilitarian and Nash group fairness, we immediately see that algorithms give dramatic improvements in 
fairness over the KPD and myopic algorithms.
All the learning algorithms also improve egalitarian group fairness over the myopic algorithm. 
But that is not the case for the KPD algorithm.
The reason for this is the KPD point system heavily prioritizes hard-to-match 
patients (recall Table~\ref{tab:kpd}), that is, it is implicitly designed to provide high egalitarian group fairness. Despite this,  
four of the learning algorithms provide better egalitarian group fairness than the KPD algorithm, 
with {\tt Exp(3)} proffering the best improvement with a fairness measure of~$1.17$.

Another valuable observation is that KPD$^+$ Pareto dominates KPD
in terms of shorter queue sizes for every group (see Table~\ref{tab:queues}). In fact, not only does each group improve but each improves 
by a large amount, with the worst case improvement being 18\%.
This algorithm then produces the best fairness of $1.22$ with respect to our primary measure, egalitarian group fairness.
Moreover, because KPD$^+$ uses the same weights for patient-donor nodes as KPD, the distribution of
groups within the waiting lists are similar. This explains its similar fairness scores of $1.28$ and $1.27$ for
 utilitarian and Nash group fairness, respectively. In contrast, our learning algorithms lead to quite different
distributions of groups within the waiting lists than the KPD algorithm; this leads to a greater fluctuations
in their three group fairness measures.
Four of the learning algorithms provide improved fairness over the KPD algorithm for all three measures.
Thus, we obtain our next lessons: \\[-0.3cm]

\noindent {\tt Lesson V.} {\em Fair algorithms can be used to reduce the accumulation of hard-to-match patients in waiting pools.}\\[-0.3cm]

\noindent {\tt Lesson VI.} {\em Fair algorithms can be implemented without detrimental affects on utilitarian measures such as waiting times and the number of transplants.}

\section{Conclusion}

In this work, we proposed a kidney exchange algorithm based on learned weights. These weights were determined via a learning approach driven by group fairness defined in terms of blood type and cPRA rates. Then we presented results obtained from simulations of dynamic KEPs, comparing the KPD algorithm used in Canada, the myopic algorithm, and our learning algorithms. Our results provide lessons for increasing  the number of transplants, decreasing waiting times, and improving group fairness. We hope this work helps contribute to enhance equity in transplantation access.

Many avenues to explore remain. First, we believe our learning algorithms can be improved. For example, the update rules for weights are rather ad-hoc. Do more structured approaches
lead to better outcomes? For example, we used the weighted power mean to evaluate the fairness of our algorithms. So can optimal weights be learnt by optimizing the weighted power mean directly during updates?

Our algorithms were tested using realistic but simulated data. The next step is to confirm our results with real patient data.
This would also allow for verification that the incorporation of omitted characteristics such as
departure rates, geographic location, age, etc, do not affect our conclusions. It is also important to confirm our lessons apply to other kidney exchange programs based on weighting systems and whether or not they extend to programs based upon hierarchical optimization. 

Improvements may also be possible by changing the set of objects that are weighted. 
Our algorithms learn node weights. Can better results be obtained by learning arc weights? That is, if a specific value 
(positive or negative) is given for the use of each specific donor-patient transplant.
In this paper, we assigned a fixed (negative) weight for each altruist. This weight was obtained via experiments over a range
of weights. This approach was vital in understanding the
importance of altruistic donors and how best to use them. But can the most effective negative altruist weight be found applying the learning method used for the positive patient-donor node weights?
Furthermore, would it be helpful to incorporate a more-refined weighting
system for the altruists? For example, it would be natural for an altruistic donor of blood group O to have a 
higher price.

Of course, learning algorithms can be used for objectives other than fairness, or for a mix of objectives.
Can dynamic learning algorithms be used to improve performance for other objective functions?
Finally, the model can be used to test the effectiveness of other proposed modifications to a KEP.
For example, currently if a donor-patient pair is compatible then that transplant is scheduled automatically.
Thus, such pairs are not included in the KPD program. Theoretically, the inclusion of
compatible donor-patient pairs may increase the chances of matching other patients in the program~\cite{SUY20}. In particular, if the patient is easy-to-match then the corresponding donor may induce analogous benefits to that of an
altruistic donor. In addition, as discussed, the donor at the end of a path is currently redesignated as a living donor in the deceased donor program. Quantifying the potential benefits of this donor remaining as a bridge donor in the KPD program is important.

\noindent {\bf Acknowledgements.} The authors thank William Klement and Mike Gillissie of Canadian Blood Services for numerous discussions and expert advice. We are also extremely grateful to  David Manlove and John Dickerson for detailed comments and advice. This project was partially supported by the
Natural Sciences and Engineering Research Council of Canada and the Institut de valorisation des donn\'ees  and Fonds de recherche du Qu\'ebec via an FRQ-IVADO Research Chair.

\bibliographystyle{plain}
\bibliography{ARXIV/arxiv}

\begin{thebibliography}{10}

\bibitem{ABT07}
D.~Abraham, A.~Blum, and T.~Sandholm.
\newblock Clearing algorithms for barter exchange markets: Enabling nationwide
  kidney exchanges.
\newblock In {\em Proceedings of the 8th ACM Conference on Electronic Commerce
  (EC)}, pages 295--304, 2007.

\bibitem{AAG15}
R.~Anderson, D.~Ashlagi, I.~Gamarnik, and A.~Roth.
\newblock Finding long chains in kidney exchange using the traveling salesman
  problem.
\newblock {\em PNAS}, 112(3):663--668, 2015.

\bibitem{ACM61}
K.~Arrow, H.~Chenery, B.~Minhas, and R.~Solow.
\newblock Capital-labor substitution and economic efficiency.
\newblock {\em Review of Economics and Statistics}, 43(3):225--250, 1961.

\bibitem{AGRR12}
I.~Ashlagi, D.~Gamarnik, M.~Rees, and A.~Roth.
\newblock The need for (long) chains in kidney exchange.
\newblock Technical report, National Bureau of Economic Research, 2012.

\bibitem{BFT11}
D.~Bertsimas, V.~Farias, and N.~Trichakis.
\newblock The price of fairness.
\newblock {\em Operations Research}, 59(1):17--31, 2011.

\bibitem{BHA19}
P.~Biro, B.~Haase-Kromwijk, T.~Andersson, E.~Asgeirsson, T.~Baltesova,
  I.~Boletis, C.~Bolotinha, G.~Bond, G.~Bohmig, L.~Burnapp, K.~Cechlarova,
  P.~Di~Ciaccio, J.~Fronek, K.~Hadaya, A.~Hemke, C.~Jacquelinet, R.~Johnson,
  R.~Kieszek, D.~Kuypers, R.~Leishman, M-A. Macher, D.~Manlove, G.~Menoudakou,
  M.~Salonen, B.~Smeulders, V.~Sparacino, F.~Spieksma, M.~Valentin, N.~Wilson,
  and J.~Van~de Klundert.
\newblock Building kidney exchange programmes in {E}urope: {A}n overview of
  exchange practice and activities.
\newblock {\em Transplantation}, 103(7):1514--1522, 2019.

\bibitem{BKM21}
P.~Biro, J.~Van~de Klundert, D.~Manlove, W.~Pettersson, T.~Andersson,
  E.~Asgeirsson, L.~Burnapp, P.~Chromy, P.~Delgado, P.~Dworczak, B.~Haase,
  A.~Hemke, R.~Johnson, X.~Klimentova, D.~Kuypers, A.~Nanni~Costa,
  B.~Smeulders, F.~Spieksma, M.~Valentin, and J.~Van~der Klundert.
\newblock Modelling and optimisation in {E}uropean kidney exchange programmes.
\newblock {\em European Journal of Operational Research}, 291(2):447--456,
  2021.

\bibitem{CBS14}
Canadian {B}lood {S}ervices.
\newblock Donation and transplantation kidney paired donation program data
  report 2009--2013.
\newblock
  \url{https://professionaleducation.blood.ca/sites/msi/files/Canadian-Blood-Services-KPD-Program-Data-Report-2009-2013.pdf},
  2014.

\bibitem{Bul03}
P.~Bullen.
\newblock {\em Handbook of Means and their Inequalities}.
\newblock Springer, 2003.

\bibitem{CBS18}
{C}anadian~{B}lood {S}ervices.
\newblock Interprovincial organ sharing national data report: Kidney paired
  donation program 2009-2018.
\newblock
  \url{https://professionaleducation.blood.ca/sites/default/files/kpd-eng\_2018.pdf},
  2018.

\bibitem{CKK12}
I.~Caragiannis, C.~Kaklamanis, P.~Kanellopoulos, and M.~Kyropoulou.
\newblock The efficiency of fair division.
\newblock {\em Theory of Computing Systems}, 50:589--610, 2012.

\bibitem{CKM19}
I.~Caragiannis, D.~Kurokawa, H.~Moulin, A.~Procaccia, N.~Shah, and J.~Wang.
\newblock The unreasonable fairness of maximum {N}ash welfare.
\newblock {\em ACM Transactions on Economics and Computation}, 7(3):1--32,
  2019.

\bibitem{CL22}
M.~Carvalho and A.~Lodi.
\newblock A theoretical and computational equilibria analysis of a multi-player
  kidney exchange program.
\newblock {\em European J. of Operational Research}, 2022.

\bibitem{CM10}
I.-K. Cho and S.~Meyn.
\newblock Efficiency and marginal cost pricing in dynamic competitive markets
  with friction.
\newblock {\em Theoretical Economics}, 5:215--239, 2010.

\bibitem{CNC15}
E.~Cole, P.~Nickerson, P.~Campbell, K.~Yetzer, N.~Lahaie, J.~Zaltzman, and
  J.~Gill.
\newblock The {C}anadian kidney paired donation program: {A} national program
  to increase living donor transplantation.
\newblock {\em Transplantation}, 99(5):985--990, 2015.

\bibitem{CON13}
M.~Constantino, X.~Klimentova, A.~Viana, and A.~Rais.
\newblock New insights on integer-programming models for the kidney exchange
  problem.
\newblock {\em European Journal of Operational Research}, 231(1):57--68, 2013.

\bibitem{DGG20}
M.~Delorme, S.~Garc{\'\i}a, J.~Gondzio, J.~Kalcsics, D.~Manlove, and
  W.~Pettersson.
\newblock New algorithms for hierarchical optimisation in kidney exchange
  programmes.
\newblock {\em Technical report ERGO 20--005, Edinburgh Research Group in
  Optimization}, 2020.

\bibitem{DGG22}
M.~Delorme, S.~Garcia, J.~Gondzio, J.~Kalcsics, D.~Manlove, W.~Pettersson, and
  J.~Trimble.
\newblock Improved instance generation for kidney exchange programmes.
\newblock {\em Computers and Operations Research}, 141, 2022.

\bibitem{DMP16}
J.~Dickerson, D.~Manlove, P.~Plaut, T.~Sandholm, and J.~Trimble.
\newblock Position-indexed formulations for kidney exchange.
\newblock In {\em Proceedings of the 17th ACM Conference on Economics and
  Computation (EC)}, pages 25--42, 2016.

\bibitem{DPS12}
J.~Dickerson, A.~Procaccia, and T.~Sandholm.
\newblock Dynamic matching via weighted myopia with application to kidney
  exchange.
\newblock In {\em Proceedings of the 26th Conference on Artificial Intelligence
  (AAAI)}, pages 1340--1346, 2012.

\bibitem{DPS12b}
J.~Dickerson, A.~Procaccia, and T.~Sandholm.
\newblock Optimizing kidney exchange with transplant chains: Theory and
  reality.
\newblock In {\em Proceedings of the 11th International Conference on
  Autonomous Agents and Multiagent Systems (AAMAS)}, pages 711--718, 2012.

\bibitem{DPS14}
J.~Dickerson, A.~Procaccia, and T.~Sandholm.
\newblock Price of fairness in kidney exchange.
\newblock In {\em Proceedings of the 13th International Conference on
  Autonomous Agents and Multiagent Systems (AAMAS)}, pages 1013--1020, 2014.

\bibitem{DS15}
J.~Dickerson and T.~Sandholm.
\newblock Future{M}atch: {C}ombining human value judgments and machine learning
  to match in dynamic environments.
\newblock In {\em Proceedings of the 29th Conference on Artificial Intelligence
  (AAAI)}, pages 622--628, 2016.

\bibitem{DGHR18}
Y.~Ding, D.~Ge, S.~He, and C.~Ryan.
\newblock A nonasymptotic approach to analyzing kidney exchange graphs.
\newblock {\em Operations Research}, 66(4):918--935, 2018.

\bibitem{DS77}
A.~Dixit and J.~Stiglitz.
\newblock Monopolistic competition and optimum product diversity.
\newblock {\em American Economic Review}, 67(3):297--308, 1977.

\bibitem{FAB21}
G.~Farnadi, W.~St-Arnaud, B.~Babaki, and M.~Carvalho.
\newblock Individual fairness in kidney exchange programs.
\newblock In {\em Proceedings of the 35th Conference on Artificial Intelligence
  (AAAI)}, pages 11496--11505, 2021.

\bibitem{FBS20}
R.~Freedman, J.~Borg, W.~Sinnott-Armstrong, J.~Dickerson, and V.~Conitzer.
\newblock Adapting a kidney exchange algorithm to align with human values.
\newblock {\em Artificial Intelligence}, 283, 2020.

\bibitem{GKW14}
K.~Glorie, J.~Van~de Klundert, and A.~Wagelmans.
\newblock Kidney exchange with long chains: {A}n efficient pricing algorithm
  for clearing barter exchanges with branch-and-price.
\newblock {\em Manufacturing and Service Operations Management},
  16(4):498--512, 2014.

\bibitem{HKJ08}
K.~Huh, M.~Kim, M.~Ju, H.~Chang, H.~Ahn, S.~Lee, J.~Lee, S.~Kim, Y.~Kim, and
  K.~Park.
\newblock Exchange living-donor kidney transplantation: Merits and limitations.
\newblock {\em Transplantation Proceedings}, 86:430--435, 2008.

\bibitem{LAM20}
E.~Lam and V.~Mak-Hau.
\newblock Branch-and-cut-and-price for the cardinality-constrained multi-cycle
  problem in kidney exchange.
\newblock {\em Computers \& Operations Research}, 115:104852, 2020.

\bibitem{MC14}
S.~Malik and E.~Cole.
\newblock Foundations and principles of the {C}anadian living donor paired
  exchange program.
\newblock {\em Canadian Journal of Kidney Health and Disease}, 1, 2014.
\newblock {\#}6.

\bibitem{MM14}
D.~Manlove and G.~O'Malley.
\newblock Paired and altruistic kidney donation in the {UK}: {A}lgorithms and
  experimentation.
\newblock {\em ACM Journal of Experimental Algorithmics}, 19(2):663--668, 2014.
\newblock {\#}2.6.

\bibitem{MGM06}
R.~Montgomery, S.~Gentry, W.~Marks, D.~Warren, J.~Hiller, J.~Houp, A.~Zachary,
  J.~Melancon, W.~Maley, H.~Rabb, C.~Simpkins, and D.~Segev.
\newblock Domino paired kidney donation: a strategy to make best use of live
  non-directed donation.
\newblock {\em The Lancet}, 368(9533):419--421, 2006.

\bibitem{PMK99}
K.~Park, J.~Moon, S.~Kim, and Y.~Kim.
\newblock Exchange-donor program in kidney transplantation.
\newblock {\em Transplantation Proceedings}, 31(1-2):356--357, 1999.

\bibitem{PCP11}
S.~Patel, P.~Chadha, and V.~Papalois.
\newblock Expanding the live kidney donor pool: Ethical considerations
  regarding altruistic donors, paired and pooled programs.
\newblock {\em Experimental and Clinical Transplantation}, 1:181--186, 2011.

\bibitem{Rap86}
F.~Rapaport.
\newblock The case for a living emotionally related international kidney donor
  exchange registry.
\newblock {\em Transplantation Proceedings}, 18:5--9, 1986.

\bibitem{RSU04}
A.~Roth, T.~Sönmez, and U.~Ünver.
\newblock Pairwise kidney exchange.
\newblock {\em Quarterly Journal of Economics}, 119(2):457--488, 2004.

\bibitem{RSU07}
A.~Roth, T.~Sönmez, and U.~Ünver.
\newblock Efficient kidney exchange: Coincidence of wants in markets with
  compatibility-based preferences.
\newblock {\em American Economic Review}, 97(3):828--851, June 2007.

\bibitem{RSU06}
A.~Roth, T.~Sönmez, U.~Ünver, F.~Delmonico, and S.~Saidman.
\newblock Utilizing list exchange and nondirected donation through `chain'
  paired kidney donations.
\newblock {\em American Journal of Transplantation}, 6:2694--2705, 2006.

\bibitem{Sol56}
R.~Solow.
\newblock A contribution to the theory of economic growth.
\newblock {\em The Quarterly Journal of Economics}, 70(1):65--94, 1956.

\bibitem{St22}
William St-Arnaud, Margarida Carvalho, and Golnoosh Farnadi.
\newblock Adaptation, comparison and practical implementation of fairness
  schemes in kidney exchange programs, 2022.
\newblock arXiv:2207.00241.

\bibitem{SU13}
T.~Sönmez and U.~Ünver.
\newblock Market design for kidney exchange.
\newblock In N.~Vulkan, A.~Roth, and Z.~Neeman, editors, {\em The Handbook of
  Market Design}, chapter~4. Oxford University Press, 2013.

\bibitem{SUY20}
T.~Sönmez, U.~Ünver, and B.~Yenmez.
\newblock Incentivized kidney exchange.
\newblock {\em American Economic Review}, 110(7):2198--2224, 2020.

\bibitem{WDA10}
E.~Woodle, J.~Daller, M.~Aeder, R.~Shapiro, T.~Sandholm, V.~Casingal,
  D.~Goldfarb, R.~Lewis, J.~Goebel, and M.~Siegler.
\newblock Ethical considerations for participation of nondirected living donors
  in kidney exchange programs.
\newblock {\em American Journal of Transplantation}, 10:1460--1467, 2010.

\end{thebibliography}



\end{document}